\documentclass[12pt]{iopart}

\usepackage{amssymb, amsfonts}
\usepackage[all]{xy}
\usepackage{bm}

\long\def\nodo#1{{}}

\def\MR#1{} 

\begin{document}

\title{\v{C}ech cocycles for quantum principal bundles}

\author{Zoran \v{S}koda}
\address{Institute Rudjer Bo\v{s}kovi\'c, Division of Theoretical Physics, \\
Bijeni\v{c}ka 54, P.O.Box 180, HR-10002 Zagreb, Croatia}
\ead{zskoda@irb.hr}
\begin{abstract}
In other to study connections and gauge theories on noncommutative spaces
it is useful to use the local trivializations of principal bundles.
In this note we show how to use noncommutative localization theory
to describe a simple version of cocycle data for the
bundles on noncommutative schemes with Hopf algebras in the
role of the structure group which locally look like Hopf algebraic
smash products. We also show how to use
these \v{C}ech cocycles for associated vector bundles.
We sketch briefly some examples related to quantum groups.
\end{abstract}
\section{Introduction}

It is often in noncommutative geometry that the Hopf algebras
appear in the role of symmetry groups (\cite{Drinfeld,Majid}).
Thus one approach of developing the gauge theories on
noncommutative spaces is to develop first the theory of principal
bundles whose structure group is replaced by a Hopf algebra and
then to develop the concept of connection on such a principal
bundle, and finally to study the action functionals developed on
this basis (\cite{BrzMaj:quantumgauge}). In noncommutative
geometry based on $C^*$-algebras, it is sufficient to model spaces
with a single global algebra, while in the algebraic setting, a
single algebra corresponds just to an affine algebraic variety.
For a more general variety, one needs to glue local (affine)
charts (\cite{Maszczyk,Ros:NcSch,vOystWill:Groth}. One flavour of
such gluing is the gluing of categories of modules via
noncommutative localization
(\cite{Ros:NcSch,skoda:nloc,skoda:gmj,vOystWill:Groth}).

In the case of noncommutative affine varieties playing the role
both of total and base space of a principal bundle, the generally
accepted notion of a Hopf algebraic principal bundle is a
Hopf-Galois extension~\cite{Majid,Montg,Schn:prinHomSp} and its
coalgebra generalizations (\cite{BrzMajid:ent}); the trivial
bundle on the other hand is the special case of Hopf algebraic
smash product. In this affine case, there are sensible and much
studied proposals (\cite{BrzMaj:quantumgauge}) how to introduce
connections, gauge transformations, curvature etc. What is not
clear in the published literature is how to glue such data between
local charts. Related problem is how to generalize these data to
nonaffine spaces, and in particular of noncommutative
schemes~\cite{Ros:NcSch}. Thus we are interested in gluing data
for nonaffine noncommutative principal bundles first and in the
next step of gluing the sections of associated vector bundles. The
connection forms may be treated as a special kind of such
sections, in the presence of a well behaved noncommutative
differential calculus.

In an ongoing project with Gabi B\"ohm, we study a global generalization of
Hopf-Galois extensions and also the corresponding cocycle data at
the categorical level: the cocycles there are made out of functors
(\cite{BohmSkoda}). That is a clean general approach,
but it is difficult to use it directly. In this note
I will present how to work out directly a special case,
where one can define and build cocycles
at the level of homomorphisms of algebras, rather than functors.
This direct approach will be likely easier to use in physical situations.

\section{Equivalences of smash products}

The basic ingredient of our construction of cocycles comes from the analysis of
equivalences of Hopf smash products. Even more general case of equivalences
of cocycled crossed products is fully characterized by Doi (\cite{Doi}).

We use the Sweedler's notation for coproduct $\Delta(h) = \sum
h_{(1)}\otimes h_{(2)}$, and its extension to right coactions
$\rho(e) = \sum e_{(0)}\otimes e_{(1)}$ and left coactions
$\rho(v)= \sum v_{(-1)}\otimes v_{(0)}$ (\cite{Majid,Montg}).
Given an action $\triangleright$ of a bialgebra $H$ on an algebra
$U$ which is Hopf (i.e. makes it a $H$-module algebra:
$h\triangleright(u v) = \sum (h_{(1)}\triangleright u)\otimes
(h_{(2)} \triangleright v)$ and $h\triangleright 1 =
\epsilon(h)1$, \cite{Majid,Montg}), the tensor product $U\otimes
H$ (where $\otimes$ for elements is traditionally written
$\sharp$) has nontrivial multiplication $(u\sharp h)(v\sharp g) =
\sum u (h_{(1)}\triangleright v) \sharp h_{(2)}g$ and a coaction
$u\sharp h\mapsto \sum (u\sharp h_{(1)})\otimes h_{(2)}$, which
together form a structure of a right $H$-comodule algebra $U\sharp
H$ which will be referred to as the smash product. If $H$ is a
Hopf algebra with antipode $S$ (as we assume from now on) then
$\gamma\circ S$ is the inverse $\gamma^{-1}$ of $\gamma$ in the
space of linear maps $Hom(H,E)$ with respect to the convolution
product $(f_1\star f_2)(h) = \sum f_1(h_{(1)})\cdot_E
f_2(h_{(2)})$. A trivial principal $H$-bundle over $U$ will be a
left $U$-module, right $H$-comodule algebra which is isomorphic to
$U\sharp H$ as $U$-module and $H$-comodule algebra. Morphisms of
trivial bundles will preserve these structures. If $E$ is a left
$U$-module right $H$-comodule, then the existence of the
isomorphism $\xi: U\sharp H\to E$ is equivalent to the existence
of a morphism of right $H$-comodule algebras $\gamma:H\to E$,
namely $\gamma(h) = \xi(1\sharp h)$, which we call the
trivializing section. The appropriate action is then
$h\triangleright a = \sum \gamma(h_{(1)})u\gamma(Sh_{(2)})$ where
$h\in H$ and $u\in U$.

Let now $f: E_1\to E_2$ be a morphism of bundles over $U$
where $\gamma_i:H\to E_i$, $i = 1,2$ are given trivializations.
Then $f\circ \gamma_2$ is also a trivializing section of $E_2$.
Thus for comparing the trivial bundles it is enough to compare different
trivializing sections of the same bundle when $f = 1$ and $E_1 = E_2 = E$.
A generic element in $E$ can be written as $\sum_k u_k \gamma_2(h_k)$.
The element $\gamma_1(h)$ is of the form
$\sum y(h_{(1)})\gamma_2(h_{(2)})$ for some algebra map $y:H\to E$. Indeed,
$y$ can be obtained as $(1\otimes \epsilon)\xi^{-1}_2(\gamma_1(h))$ where
$\xi^{-1}:E\to U\sharp H$ is the isomorphism induced by the section $\gamma_2$.
Thus
\begin{equation}\label{eq:y}
y(h) = \sum \gamma_1(h_{(1)})\gamma_2^{-1}(h_{(2)})
\end{equation}
encodes all the information on comparing different
trivializations. Furthermore, $$\begin{array}{l}\sum
(h_{(1)}\triangleright_1 u)y(h_{(2)}) = \sum
(h_{(1)}\triangleright_1
u)\gamma_1(h_{(2)})\gamma_2^{-1}(h_{(3)})\\ = \sum
y(h_{(1)})\gamma_2(h_{(2)}) u \gamma_1(Sh_{(3)})
\gamma_1(h_{(4)})\gamma_2^{-1}(h_{(5)})= \sum
y(h_{(1)})(h_{(2)}\triangleright_2 u).\end{array}$$ Therefore
\begin{equation}\label{eq:actionsComp}
\sum (h_{(1)}\triangleright_1 u)y(h_{(2)}) =
\sum y(h_{(1)})(h_{(2)}\triangleright_2 u)
\end{equation}
Furthermore it is easy to see that a composition of
two morphisms of trivial bundles
corresponds to a convolution product of the corresponding $y$-maps:
$$
y_{13}(h) = \sum \gamma_1(h_{(1)})\gamma_3^{-1}(h_{(2)}) = \sum
\gamma_{(1)}(h_{(1)})\gamma_{(2)}^{-1}(h_{(2)})\gamma_2(h_{(3)})\gamma_3(h_{(4)})
= \sum y_{12}(h_{(1)})y_{23}(h_{(2)})
$$

\section{\v{C}ech cocycle}

\subsection{Noncommutative localization and noncommutative schemes}

We start with rather technical general sketch and then we give a
simple recipe of a special case of our interest. A noncommutative
space in algebraic sense is ultimately represented by an abelian
category $\mathcal{A}$ whose objects are viewed as quasicoherent
sheaves. The category is covered by {\it flat localization}
functors $Q_\lambda^* :\mathcal{A}\to\mathcal{A}_\lambda$,
$\lambda\in\Lambda$ having right adjoint functors $Q_{\lambda*}$,
and $\mathcal{A}_\lambda$ is isomorphic to the category of left
modules over a noncommutative algebra $U_\lambda$
(\cite{skoda:nloc,Ros:NcSch}). Denote also $Q_\lambda = Q_{\lambda
*} Q_\lambda^*:\mathcal{A}\to\mathcal{A}$. Thus we may be given a
family of algebras $U_\lambda$, $\lambda\in\Lambda$ viewed as
algebras of functions on Zariski open charts whose {\it
intersections} will be replaced by considering mixed iterates
$Q_\lambda^* Q_{\mu*}Q^*_{\mu*}:
\mathcal{A}\to\mathcal{A}_{\mu\lambda}\subset\mathcal{A}_\lambda$
of localization functors where $\mathcal{A}_{\mu\lambda}$ is the
essential image of the written functor with values in
$\mathcal{A}_\lambda$. The bad thing is that
$\mathcal{A}_{\mu\lambda}$ is not necessarily equivalent to
$\mathcal{A}_{\lambda\mu}$, nor $Q_\lambda Q_\mu\cong Q_\mu
Q_\lambda$, and worse, with $A_\lambda$ looking affine, i.e. like
the category of modules over and algebra, the iterates
$\mathcal{A}_{\lambda\mu}$ do not look like that in general, hence
we can {\it not} talk about the {\it algebra} $U_{\lambda\mu}$,
but rather only a bimodule.

\subsection{Mixed localizations and gluing bimodules}

The double consecutive localization
$\mathcal{A}_{\lambda\mu}\subset\mathcal{A}_\mu$, hence
$\mathcal{A}_{\lambda\mu}$ is always a subcategory of the category
of all left $U_\mu$-modules. One can still define a left
$U_\mu$-module $U_{\lambda\mu}:= Q_\mu(U_\lambda)$. Under mild
conditions, it is a $U_\mu$-$U_\lambda$-bimodule flat from both
sides; this bimodule is generated by the image of unit element in
$U_\lambda$ under the adjunction map $U_\mu\to U_{\lambda\mu}$.
This will be the working assumption in this article. In fact, if
the family of localizations is finite, this can be abstracted
further to the noncommutative space covers of Kontsevich-Rosenberg
(\cite{KontsRos}) which are given in terms of a coring with an
additional "structure map"; it seems that the whole theory of this
article could be generalized to their setup (see
\cite{skoda:hgcoring} for some hints) In the simplest affine case,
when there is a global algebra $U$ such that $U_\lambda =
Q_\lambda(U)$ for all $\lambda$, the localization functor $Q_\mu$
is isomorphic to $U_\mu\otimes_U$ and $U_{\lambda\mu} =
U_\mu\otimes_U U_\lambda$. Notice also that we can consider higher
iterates like $U_{\lambda\mu\nu} = U_\nu\otimes_U U_{\mu}\otimes_U
U_\lambda$ and various natural maps from lower into higher
iterates. We are more interested in the case when $U_\lambda =
E_\lambda^{\mathrm{co}H}$ is the subalgebra of coinvariants in a
right $H$-comodule algebra $E_\lambda$, i.e. the elements $u\in
E_\lambda$ for which the coaction is of the form $\rho(u) =
u\otimes 1$; and such that $E_\lambda= Q_\lambda(E)$ is the
localized algebra of some algebra $E$. We say that $U_\lambda$ is
the algebra of localized $H$-coinvariants in $E$. In good cases
such algebras can be glued to form a noncommutative {\it quotient
space} of $E$, the method which we pioneered in \cite{skoda:ban}.
Now we can define $U_{\lambda\mu} = E_{\lambda\mu}^{\mathrm{co}H}$
and similarly for higher iterates. This will be important for
\v{C}ech cocycles.

\subsection{Locally trivial bundles}

Suppose we are now given system of algebras
$\{U_\lambda\}_{\lambda\in\Lambda}$ as above. We should
technically require that their categories of modules glue
appropriately along $U_{\lambda\mu}$ and $U_{\mu\lambda}$ via flat
descent to the global category $\mathcal{A}$. As this is treated
at length elsewhere we just suppose that the standard setup, e.g.
of noncommutative schemes is understood; we also follow the
assumptions from the previous section about double iterates. While
one can define a rather abstract notion of principal $H$-bundle on
$\mathcal{A}$ (\cite{BohmSkoda}) we will not do that here (but the
comparison is due in a sequel paper). Instead we propose the new
explicit concept of {\bf noncommutative \v{C}ech 1-cocycle with
coefficients in the Hopf algebra} $H$:

(1) Each $U_\lambda$ is equipped with a Hopf action
$H\triangleright_\lambda U\to U$

(2) For each ordered pair $\lambda,\mu$ there is a map
$y_{\mu\lambda}: H\to U_{\lambda\mu}$
(notice the order of labels and that the codomain is
not generally an algebra but a bimodule) such that

\begin{equation}
\label{eq:actionsCoc}
\sum (h_{(1)}\triangleright_\mu u)y_{\mu\lambda}(h_{(2)}) =
\sum y_{\mu\lambda}(h_{(1)})(h_{(2)}\triangleright_\lambda u)
\end{equation}
Notice that the multiplications on the left and right
are due the bimodule structure on $U_{\lambda\mu}$.

(3) For each ordered triple $\lambda,\mu,\nu$ we have the cocycle conditions:
for each $h\in H$,
$$
y_{\nu\lambda}(h) = \sum y_{\nu\mu}(h_{(1)})y_{\mu\lambda}(h_{(2)})
$$
in $U_\nu$-$U_\lambda$-bimodule $U_{\lambda\mu\nu}$ and
$$
y_{\mu\mu}(h) = 1
$$
in $U_{\mu\mu}\cong U_\mu$.

Consider another set of such data with
the same cover $\{U_\lambda\}_{\lambda\in\Lambda}$ (gluing data for
the base space understood), but different $H$-actions
$\tilde\triangleright_\lambda$ and different bundle transition
maps $\tilde{y}_{\lambda\mu}$. The two sets are cohomologous if there is
a 0-cocycle which is a family of linear maps $r_\lambda:H\to U_\lambda$
and this 0-cocycle relates the above 1-cocycles as follows:

(1) each $r_\lambda$ is convolution invertible

(2) $r_\lambda$ exhibits the equivalence of $\triangleright_\lambda$
actions $\triangleright_\lambda$, $\tilde\triangleright_\lambda$:
\begin{equation}\label{eq:cohom1}
\sum (h_{(1)}\triangleright_\lambda u)r_\lambda(h_{(2)}) =
\sum r_\lambda(h_{(1)})(h_{(2)}\tilde\triangleright_\lambda u)
\end{equation}

(3) for each ordered pair $(\mu,\lambda)$,
\begin{equation}\label{eq:cohom2}
\sum r_\lambda(h_{(1)}) y_{\lambda\mu}(h_{(2)})
= \sum y_{\lambda\mu}(h_{(1)})r_\mu(h_{(2)})
\end{equation}
holds in $U_{\mu\lambda}$.

\section{Examples over quantum flag varieties}

Given an indeterminate $q$, $M_q(2)$ is the associative
$\mathbf{C}[q,q^{-1}]$-algebra with generators $a=
t^1_1,b=t^1_2,c=t^2_1,d=t^2_2$ subject to the relations $ab =
qba$, $ac = qca$, $bd = qdb$, $cd = qdc$, $bc = cb$ and $ad - da =
(q-q^{-1})bc$. It generalizes to $M_q(n)$ which is generated by
$t^i_j$ where $i,j= 1,\ldots,n$ and for every pair $i< j$, $k< l$
of labels modulo any set of relations such that setting $a=t^i_k$,
$b=t^i_l$, $c=t^j_k$, $d = t^j_l$ generates the subalgebra which
is a copy of $M_q(2)$. It is convenient to form the matrix $T =
(t^i_j)$. The algebra $GL_q(n)$ is the localization fo $M_q(n)$ at
the central element, the quantum determinant $D =
\sum_{\sigma\in\Sigma} (-q)^{l(\sigma)} t^1_{\sigma(1)}
t^2_{\sigma(2)}\cdots t^n_{\sigma(n)}\in M_q(n)$ where $l$ is the
length of the permutation $\sigma$. Similarly one defines quantum
minors for the submatrices of $T$.

The formulas $\Delta(t^i_j) =\sum_{k=1}^n t^i_k\otimes t^k_j$,
$\epsilon(t^i_j) = \delta^i_j$ uniquely extend to homomorphisms of
algebras making $GL_q(n)$ a Hopf algebra with an antipode $S$ such
that $S T = T^{-1}$. The subalgebra generated by all $t^i_j$ with
$i<j$, is a Hopf ideal $I$ and the quotient Hopf algebra will be
referred to as quantum Borel subgroup $B_q$ with generators
$h^i_j=t^i_j+I$ and the projection map $\pi:GL_q(n)\to B_q(n)$ is
given by $t^i_j\mapsto h^i_j$. $B_q$-coaction $\rho =
(\mathrm{id}\otimes\pi)\circ\Delta:GL_q(n)\to GL_q(n)\otimes B_q$
makes $GL_q(n)$ a right $B_q$-comodule algebra. We have earlier
exhibited a family (\cite{skoda:ban}) of $n!$ Ore localizations of
$GL_q(n)$ (Ore condition follows from \cite{skoda:ore}) which are
$\rho$-compatible (\cite{skoda:ban,skoda:gmj}) in the sense that
the coaction extends to the localized algebra in a unique way
making it a $B_q$-comodule algebra.

The construction of a coset space symbolically denoted by
$GL_q(n)/B_q(n)$ (or its $SL_q(n)$-version) is coming in a packet
with the local trivialization: both stem from the geometric
understanding of the quantum Gauss decomposition
(\cite{skoda:ban,skoda:qGauss}). Namely, one decomposes the matrix
of the generators $T$ with rows permuted by the permutation matrix
$w_{\sigma}^{-1} T$ as $U_\sigma A_\sigma$ where $U_\sigma$ is an
upper diagonal unidiagonal matrix and $A_\sigma$ the lower
triangular matrix, both with entries in the quotient skewfield of
$GL_q(n)$. Here $\sigma\in\Sigma(n)$ is the element of the
permutation group (the Weyl group for our case) and there are $n!$
such elements. The entries of $U_\sigma$ and of $A_\sigma$
together generate a subalgebra of the quotient field which is
isomorphic to the Ore localization of $G = GL_q(n)$ by the set of
principal (=lower right corner) quantum minors of $w_{\sigma}^{-1}
T$. This localization $G_\sigma$ is compatible with the coaction
of quantum Borel so we have the induced coaction $\rho_\sigma$
which makes $G_\sigma$ a right $B_q(n)$-comodule algebra. The
entries of $U_\sigma$ are coinvariant hence form a chart on the
quantum homogeneous space (for more precise statement see
\cite{skoda:ban,skoda:qGauss}). On the other hand,
$\gamma_\sigma:B_q(n)\to G_\sigma$ is defined by the simple rule
on generators
$$
\gamma_\sigma(h^i_j) = (A_\sigma)^i_j
$$
and extended as a homomorphism of algebras. To find the
$\gamma_\sigma$ hence it suffices to know how to do the Gauss
decomposition of matrices with noncommutative entries and for the
convolution inverse, one needs in addition to compute the
antipode. Both problems can be done in terms of quasideterminants
of Gel'fand and Retakh (\cite{GelfandRetakh,skoda:nloc}), and the
more detailed formulas are left for \cite{skoda:qGauss}. The
quasideterminant involved, that is the quasiminors of $T$ with
some rows permuted, can further be expressed (up to a factor of
$(-q)$ to some power) as a ratio of two quantum determinants.

The Hopf action on the coinvariants is given on generators by
$$
h^i_j\triangleright_\sigma (U_\sigma)^k_l = \sum_{i\geq s\geq j}
\gamma_\sigma(h^i_s)(U_\sigma)^k_l \gamma_\sigma(S h^s_j)
$$
where $\gamma_\sigma$ is a ratio of quasideterminants
$$
\gamma_\sigma(h^i_s) =
|T^{\sigma(i),\sigma(j+1),\sigma(j+2),\ldots,\sigma(n)}_{j,j+1,j+2,\ldots,j+n}
|_{ij}|T^{\sigma(j),\sigma(j+1),\sigma(j+2),\ldots,\sigma(n)}_{j,j+1,j+2,\ldots,j+n}|^{-1}_{jj}
$$
Thus we can easily find $y_{\sigma\tau}(h) = \sum
\gamma_\sigma(h_{1)})\gamma_\tau(Sh_{(2)})$ on the generators.
Notice that $y_{\sigma\tau}$ is not the homomorphism of algebras
so one needs to go back to $\gamma$-s to compute it on an an
arbitrary given element.

For example on the simplest case of $G = SL_q(2)$ one has
$$
\gamma_{I}(h^i_j) =
\left(\begin{array}{cc}a-bd^{-1}c&0\\c&d\end{array}\right)
,\,\,\,\,\,\,\,\gamma_{\tau}(h^i_j)
=\left(\begin{array}{cc}c-db^{-1}a&0\\a&b\end{array}\right)
$$
where $I$ is the trivial and $\tau$ the nontrivial permutation on
two letters, while the matrix of transition maps is
$$
Y = \left(\begin{array}{cc}y(h^1_1) & y(h^1_2)\\
y(h^2_1)& y(h^2_2)\end{array}\right)
= \left(\begin{array}{cc}-u & 0 \\
1 & u'\end{array}\right),\,\,\,\,\,\,Y^{-1} = \left(\begin{array}{cc} -u'&0 \\
u & 0 \end{array}\right),
$$
where $u = bd^{-1}$ is the generator in chart of the homogeneous
space corresponding to the trivial permutation and $u'= db^{-1}$,
the generators of the chart corresponding to the nontrivial
permutation in $\Sigma(2)$. Though the base looks like
$\mathbf{C}P^1$ at the local algebra level, its further structures
are nonclassical (e.g. the measure utilized in \cite{skoda:cs}).

\section{The associated vector bundles}

For a right $H$-comodule algebra $E$
and a left $H$-comodule $V$ with coactions
$\rho_E:E\to E\otimes H$ and $\rho_V: V\to H\otimes V$,
the cotensor product
is the vector subspace $E\Box V\subset E\otimes V$
which equalizes $\rho_E\otimes \mathrm{id}_V$ and
$\mathrm{id}_E\otimes \rho_V$.
If $E$ is a faithfully flat Hopf-Galois extension of $E^{\mathrm{co}H}$ then
$E$ may be interpreted as a principal bundle and $E\quad V$
as the space of global sections of
the associated vector bundle with typical fiber $V$.
We sketched in \cite{skoda:cs} that the cotensor products may be glued as
well. The transition cocycles from the previous section may be used as
well.

Let $v_j$, $j=1,\ldots,n$ be a basis of $V$ (for simplicity,
we consider the finite-dimensional fiber) and the coaction $\rho_V$
in this basis is given by
$$
\rho_V(v_i) = \sum_j v^j_i\otimes v_j
$$
for some $v^j_i\in H$.
The coaction axiom implies that
$\Delta(v^j_i)=\sum_k v^k_i\otimes v^j_k$
(notice that the matrix multiplication is transposed).
Consider now a cocycle for the principal bundle
with the notation from above. Then
$$
\gamma_\lambda(v^j_i) = \sum_k y_{\lambda\mu}(v^k_i)\gamma_\mu(v^j_k)
$$
Define the transition matrices $\mathfrak{M}_{\lambda\mu}$ by
$$
(\mathfrak{M}_{\lambda\mu})^i_j = y_{\lambda\mu}(v^i_j)\in U_{\mu\lambda}.
$$
Then $\mathfrak{M}_{\lambda\mu}^T\mathfrak{M}^T_{\mu\nu} =
\mathfrak{M}^T_{\lambda\nu}$ and
$\mathfrak{M}_{\lambda\lambda}=I$, where one interprets the
results in appropriate iterated localizations and where $()^T$ is
the sign for transposition (in $i\leftrightarrow j$).

Triviality of a Hopf-Galois extension $E^{\mathrm{co}H}\hookrightarrow E$
implies triviality of the associated bundle in the following sense.
If a right $H$-comodule $E$ admits convolution invertible map of
$H$-comodules $\gamma:H\to E$
then for any left $H$-comodule $V$ there is
an automorphism of left $E^{\mathrm{co}H}$-modules
$$
\kappa^\gamma_V : E\otimes V\cong E\otimes V,\,\,\,\,
\kappa^\gamma_V(e\otimes v) = \sum e\gamma(v_{(-1)})\otimes v_{(0)}
$$
with inverse $\bar\kappa^\gamma_V:
e\otimes v\mapsto \sum e\gamma^{-1}(v_{(-1)})\otimes v_{(0)}$,
and the automorphism $\kappa^\gamma_V$ restricts to
the isomorphism of $E^{\mathrm{co}H}$-modules
$$
\kappa^\gamma_V|:E^{\mathrm{co}H}\otimes V\to E\Box V.
$$
Now if the data $U_\lambda$, $\triangleright_\lambda$, $y_{\lambda\mu}$
form a cocycle of a principal $H$-bundle, then we can define the space of
global sections of the associate bundle with typical fiber $V$
as the vector subspace $\Gamma \xi_V$ of
$$
\prod_{\lambda\in\Lambda} U_\lambda\otimes V
$$
consisting of $|\Lambda|$-tuples $(\sum_{i = 1}^{n_\lambda}
u_i^\lambda\otimes v_i)_\lambda$ such that $\sum_i u_i^\mu
y_{\mu\lambda}(v_{i(-1)})\otimes v_{i(0)} = \sum_j
u_j^\lambda\otimes v_{j}$ in $U_{\lambda\mu}\otimes V$ for all
ordered pairs $(\lambda,\mu)$.

\section{Conclusion: perspective toward gluing connections}

As we know how to define the global sections of associated vector
bundles, in particular we can do that for defining connections
globally from local pieces. Available definitions of connections
(see e.g. \cite{BrzMaj:quantumgauge}) use as an ingredient the
differential calculi over noncommutative algebras. The question of
gluing the noncommutative calculi itself has some new elements,
for example there is an additional condition of compatibility of a
differential calculus with localizations involved in a cover. For
that reason I have introduced in \cite{skoda:nloc} the notion of
{\it differential Ore condition} which might be useful.

Some nonaffine examples could already be constructed for
homogeneous spaces of quantum groups. Some natural differential
calculi there are known \cite{Jurco:calc,Wor:calc} and we saw in
this article how to define local trivializations for certain
canonical principal bundles over quantum flag variety.



\section*{References}
\begin{thebibliography}{10}
\bibitem{BohmSkoda}
G. B\"ohm, Z. \v{S}koda, Globalizing Hopf-Galois extensions, in preparation

\bibitem{BrzMaj:quantumgauge}
T.~Brzezi\'nski, S.~Majid, Quantum group gauge theory
on quantum spaces, Comm. Math. Phys. 157 (1993), no. 3, 591--638.

\bibitem{BrzMajid:ent}
T.~Brzezi\'nski, S.~Majid, Coalgebra bundles, Comm.
Math. Phys.~191 (1998), no. 2, 467--492.

\bibitem{Doi}
Y. Doi, Equivalent crossed products for a Hopf
algebra, Comm. Alg. {\bf 17} (1989), 3053--3085.

\bibitem{Drinfeld}
V. Drinfel'd, Quantum groups, Proc. ICM 1985

\bibitem{GelfandRetakh}
I.M. Gel'fand, V.S. Retakh, Determinants of matrices over
noncommutative rings, Funct.Anal.Appl. 21 (1991), 51–-58;
Quasideterminants I, Selecta Mathematica, N.S. 3 (1997) 4,
517–-546, {\tt q-alg/9705026}.

\bibitem{Jurco:calc}
B. Jur\v{c}o, Differential calculus on quantum groups:
Constructive procedure. in: Proc. Enrico Fermi School on Quantum
Groups, Como 1994:193--202, {\tt hep-th/9408179}.

\bibitem{KontsRos} M.~Kontsevich, A.~L. Rosenberg,
Noncommutative smooth spaces, The Gelfand Math. Seminars,
1996--1999, Birkh\"{a}user 2000, 85--108; {\tt
arXiv:math.AG/9812158}.

\bibitem{Majid} S.~Majid, Foundations of
quantum group theory, Cambridge University Press 1995.

\bibitem{Maszczyk}
T. Maszczyk, Noncommutative geometry through monoidal
categories, {\tt math.QA/0611806}

\bibitem{Montg}
S. Montgomery, Hopf algebras and their actions on
rings, CBMS Regional Conference Series in Mathematics {\bf 82},
AMS 1993.


\bibitem{Ros:NcSch} A.~L. Rosenberg,
Noncommutative schemes, Comp. Math. 112 (1998), pp.~93--125.

\bibitem{Schn:prinHomSp} H.-J.~Schneider,
Principal homogeneous spaces for arbitrary Hopf algebras, Israel
J.~Math. 72, (1990), nos.1-2, pp. 167--195.

\bibitem{skoda:ore}
Z. \v{S}koda, Every quantum minor generates an Ore set , Int.
Math. Res. Notices 2008, rnn063-8; {\tt math.QA/0604610}.

\bibitem{skoda:hgcoring}
{\tt
http://ncatlab.org/zoranskoda/show/cleft+extension+of+a+space+cover}

\bibitem{skoda:cs} Z.~\v{S}koda, Coherent
states for Hopf algebras, Letters in Mathematical Physics {\bf
81}, 1, (2007) 1--17, {\tt math.QA/0303357}.

\bibitem{skoda:qGauss} Z. \v{S}koda, Geometric meaning of quantum
Gauss decomposition, in preparation.

\bibitem{skoda:ban} Z. \v{S}koda, Localizations for
construction of quantum coset spaces, in "Noncommutative geometry
and Quantum groups", Banach Center Publications {\bf 61} (2003),
265--298, {\tt math.QA/0301090}.

\bibitem{skoda:nloc} Z. \v{S}koda,
Noncommutative localization in noncommutative geometry,  London
Math. Soc. Lecture Note Series {\bf 330}, ed. A. Ranicki; pp.
220--313, Cambridge Univ. Press, 2006, {\tt math.QA/0403276}.

\bibitem{skoda:bun} Z. \v{S}koda, Quantum bundles using coactions
and localizations, in preparation

\bibitem{skoda:gmj} Z. \v{S}koda, Some equivariant constructions
in noncommutative algebraic geometry, Georgian Math. J. {\bf 16}
(2009), n. 1, 183-202, {\tt arXiv:0811.4770}.

\bibitem{vOystWill:Groth}
F.~van Oystaeyen, L.~Willaert, Grothendieck topology, coherent
sheaves and Serre's theorem for schematic algebras, J.~Pure Appl.
Alg.~104 (1995), pp.~109--122; Cohomology of schematic algebras,
J.~Alg. 185 (1996), pp.~74--84.

\bibitem{Wor:calc}
S. L. Woronowicz, Differential calculus on compact matrix
pseudogroups (quantum groups), Commun. Math. Phys. {\bf 122},
125–-170 (1989)
\end{thebibliography}
\end{document}